\documentclass[11pt]{article}
\usepackage{amsfonts}
\usepackage{color}
\usepackage{graphics}
\usepackage{indentfirst}

\usepackage{cite}
\usepackage{latexsym}
\usepackage{amsmath}
\usepackage{amssymb}
\usepackage[dvips]{epsfig}
\usepackage{amscd}

 \usepackage{color}
 \usepackage{indentfirst}
\usepackage{ntheorem}
\theoremstyle{plain}
\theoremseparator{.}

\newtheorem{theorem}{Theorem}[section]
\newtheorem{remark}{Remark}[section]

\newtheorem{lemma}[theorem]{Lemma}

\newtheorem{corollary}[theorem]{Corollary}

\newcommand\thmref[1]{Theorem~\ref{#1}}

\newcommand{\ti}{\tilde}

\newcommand\tep{ (\te^{-1/2}-\al_1^{-1/2})_+}
\def\pf{{\it Proof.}  }

\newcommand{\thatsall}{\hfill$\Box$}
\newcommand{\bi}{\bibitem}
\newcommand{\dis}{\displaystyle}
\newcommand{\bt}{\begin{theorem}}
\newcommand{\bl}{\begin{lemma}}
\newcommand{\el}{\end{lemma}}
\newcommand{\et}{\end{theorem}}

\renewcommand{\b}{\beta  }

\newcommand{\te}{\theta}

\newcommand{\al}{\alpha}

\newcommand{\ve}{\varepsilon}
\newcommand{\la}{\label}

\newcommand{\ka}{\kappa}

\newcommand{\bn}{\begin{eqnarray}}
\newcommand{\en}{\end{eqnarray}}
\newcommand{\bnn}{\begin{eqnarray*}}
\newcommand{\enn}{\end{eqnarray*}}

\newcommand{\bnnn}{\begin{eqnarray*}}
\newcommand{\ennn}{\end{eqnarray*}}
\newcommand{\ben}{\begin{enumerate}}
\newcommand{\een}{\end{enumerate}}
\newcommand{\ba}{\begin{aligned}}
\newcommand{\ea}{\end{aligned}}
\newcommand{\be}{\begin{equation}}
\newcommand{\ee}{\end{equation}}

\def\norm[#1]#2{\|#2\|_{#1}}

\makeatletter
\@addtoreset{equation}{section}
\makeatother

\makeatletter
\@addtoreset{equation}{section}
\makeatother

\title{ Nonlinearly stability of solutions  on the outer Pressure Problem of Compressible Navier-Stokes System  with
Temperature-Dependent  Heat Conductivity
\thanks{ Partially supported by NNSFC Grant Nos. 11501143, 12071200, 11971401, 11871410 and 11871408.}}

\date{ }
\author{Guocai Cai$^a$, Canpei Chen$^b$,  Yanfang Peng$^c$\thanks{
  Email addresses: gotry@xmu.edu.cn (G.C. Cai), 20150008029@m.scnu.edu.cn (C.P. Chen), pyfang2005@sina.com (Y.F. Peng)} \\[3mm]  a. School of Mathematical Sciences, Xiamen University,\\
Xiamen 361005, P. R. China; \\b. School of Mathematical Sciences,\\
  University of Chinese Academy of Sciences, Beijing 100049, P. R. China\\ c. School of Mathematical Sciences,   Guizhou Normal University,  \\Guiyang, 550001, P. R. China; }
\begin{document}

\maketitle
\begin{abstract}
 In this paper, the one-dimensional compressible Navier-Stokes system with outer pressure boundary  conditions is investigated. Under some suitable assumptions, we prove that the specific volume and the temperature are bounded from below and above independently of time, and then give the local and global existence of strong solutions. Furthermore, we also obtain the convergence of the global strong solution to a stationary state and the nonlinearly stability  of its convergence. It is worth noticing that all the assumptions imposed on the initial data are the same as Takeyuki Nagasawa [Japan.J.Appl.Math.(1988)]. Therefore, our work can be regarded as an improvement of the results of Takeyuki.
 \end{abstract}



\noindent{\bf Keywords:}  Compressible Navier-Stokes system;
Temperature-dependent heat conductivity;  Outer pressure problem; Stability

\noindent \textbf{Math Subject Classification}: 35Q35; 76N10.

\section{Introduction}

The one-dimensional  compressible Navier-Stokes system   describing   the   motion of a
viscous heat-conducting   perfect polytropic gas has the
form in the Lagrangian mass coordinates (see \cite{4,23}):
  \be\la{1.1}
v_t=u_{x},
\ee
\be\la{1.2}
u_{t}+P_{x}=\left(\mu\frac{u_{x}}{v}\right)_{x},\ee
 \be\la{1.3}\left(e+\frac{1}{2}u^2\right)_t+\left(Pu\right)_x=\left(\frac{\kappa\te_x+\mu uu_x}{v}\right)_{x}, \ee
where $t>0$ is time, $x\in \Omega= (0,1)$ denotes the
Lagrange mass coordinate,  and the unknown functions $v>0, u$ and $P$ are the specific volume of the gas, fluid velocity and  pressure respectively. In this paper, we focus on ideal polytropic gas, that is, $P$ and $e$ satisfy
\be \la{1.4}   P =R \theta/{v},\quad e=c_v\theta +\mbox{const},
\ee
where  both specific gas constant  $R$ and   heat
capacity at constant volume $c_v $ are   positive constants. The viscosity $\mu$ and the heat conductivity $\ka$  are given by  \be\la{1.5}   0<\mu=\tilde{\mu}\theta^\alpha,\quad \ka=\ti\ka \te^\beta, \ee
with constants $  \tilde{\mu}, \ti\ka>0$ and $  \alpha, \beta\ge 0.$
The system \eqref{1.1}-\eqref{1.5} is supplemented
with the initial  data for any $x\in (0,1)$,
\be\la{1.6}(v,u,\te)(x,0)=(v_0,u_0,\te_0)(x),   \ee
and   the boundary conditions
\be\la{1.8}
\left(\frac{\mu}{v}u_x-P\right)(0,t)=\left(\frac{\mu}{v}u_x-P\right)(1,t) \triangleq -p(t)<0, \quad\te_x(0,t)=\te_x(1,t)=0.
\ee
where the first boundary conditions are called outer pressure boundary conditions.
 In addition, the initial data \eqref{1.6} should be compatible with the boundary conditions \eqref{1.8}.

A large number of literatures have been devoted to the large-time existence and behavior of solutions to one-dimensional compressible Navier-Stokes equations. For constant coefficients $(\al=\beta=0)$  with large initial data,  Kazhikhov and Shelukhin\cite{9} first showed  the global existence of solutions  under the following boundary conditions:
\be\ba \la{1.9} u(0,t) =u(1,t)=0,\quad\te_x(0,t)=\te_x(1,t)=0,\ea\ee
which mean that the gas is confined into a fixed tube with impermeable gas. From then on,
significant progress has been made on the mathematical aspect of the
initial
boundary value problems, see \cite{1,2,3,15,17,18} and the references therein.
Motivated by the fact that in the case of isentropic flow a temperature dependence on the viscosity translates
into a density dependence, there is a lot of research work (see \cite{hs1,hsj1,HuangShi,3,ssyy11,28} and the references therein) concentrating on  the case that $\mu$  is
independent of $\te$, and heat conductivity is allowed to depend on temperature in a special
way with a positive lower bound and balanced with corresponding constitution relations.
In the case that $\alpha=0, \beta>0$ with  boundary condition of either \eqref{1.8} with $p(t)\equiv 0$ or
\eqref{1.9}, Kawohl \cite{ssyy11}, Jiang \cite{ssyy22,ssyy33}  and Wang \cite{ssyy44} established the
global existence of smooth solutions for \eqref{1.1}--\eqref{1.3}, \eqref{1.6} under the assumption $\mu(v) \geq \mu_0> 0$ for any $v>0$ and $\kappa$ may depend on both density and
temperature. Moreover, when $\al=0$ and $\beta\in(0,3/2),$ Jenssen-Karper \cite{24}  proved the global existence of a weak solution to \eqref{1.1}--\eqref{1.8}.
 Later on, for $\al=0$ and $\beta\in(0,\infty),$ Pan-Zhang\cite{28} obtained  the global strong solutions under the  non-slip and heat insulated boundary conditions. Recently, Duan-Guo-Zhu \cite{7}  obtained  the global strong solutions of \eqref{1.1}-\eqref{1.8}  under the condition that
  \be \la{h1.8}( v_0,u_0,\te_0)\in H^1\times H^2\times H^2. \ee
And then,  Huang-Shi-Sun \cite{huang} relaxed the initial condition \eqref{h1.8} to
  \be \la{h1.9}( v_0,u_0,\te_0)\in H^1\times H^1\times H^1.\ee

  On the other hand, concerning the large-time behavior of the strong solutions to  \eqref{1.1}-\eqref{1.6} and \eqref{1.9}, Kazhikhov \cite{Kaz}
  first proved that for the case $\alpha=\beta=0$, the strong solution is nonlinearly exponentially stable  as time tends to
  infinity. More recently, Huang-Shi \cite{HuangShi} extended  the same result  to  $\alpha=0,\beta>0$.  However, it should be pointed out  here that the method used in \cite{Kaz}  relies heavily on the non-degeneracy of both the viscosity $\mu$ and the heat conductivity $\ka$, so that it can not be applied directly to  the degenerate and nonlinear case  ($\al=0,\b>0$).

  When it comes to the outer pressure problem \eqref{1.1}-\eqref{1.8} with $\alpha=\beta=0$, Takeyuki Nagasawa \cite{tak} showed  the convergence of the solutions to a stationary state and the rate of its convergence under some suitable assumptions. However, whether the same results are still hold  for the system \eqref{1.1}-\eqref{1.8} with $\alpha=0,\beta>0$ remains open.
One of our purposes  is to give a positive answer to this problem. In fact, we establish the existence of the strong solutions for the problem \eqref{1.1}-\eqref{1.8} and study the large-time stability under the additional assumptions \eqref{eq1.9} and \eqref{eq1.799}, more precisely, we have main results as follows.

 \begin{theorem}\la{thm1.2}
 Suppose that
 \bnn    \alpha=0, \quad  \beta> 0 ,\quad 0<p(t)\in C^1[0,+\infty) ,\enn
 and the initial data $ ( v_0,u_0,\te_0)$ satisfies
  \be\la{h1.100} ( v_0,u_0,\te_0)\in   H^1 (0,1),  \inf_{x\in [0,1]}v_0(x)>0, \inf_{x\in [0,1]}\theta_0(x)>0.\ee
Then  the initial-boundary-value problem \eqref{1.1}-\eqref{1.8} has a unique global strong solution $(v,u,\te)$ such that for each $T>0$,
 \be
 \begin{cases} \la{1.111}  v ,\,u,\,\theta \in L^\infty(0,T;H^1(0,1)),\\ v_t\in
  L^\infty(0,T;L^2(0,1))\cap L^2(0,T;H^1(0,1)), \\ u_x,\,\theta_x,\,u_t\,,\theta_t,\,v_{xt},\,\,u_{xx},\,\te_{xx} \in
  L^2((0,T)\times(0,1)).\end{cases}
  \ee
 Moreover,
  \be\la{1.13}  C^{-1}\leq v(x,t)\leq C ,\quad  C^{-1}\leq \theta(x,t)\leq C,\ee
 where $C$ is a positive constant depending on the initial data and $T$.
 \end{theorem}

The following theorem indicates that the conclusion of Theorem \ref{thm1.2} is uniform with respect to $t$  and the global strong solution tends to a stationary state as time goes away, if one further assumes that $p(t) $ satisfies
  \be\la{eq1.9}
  m_p\triangleq \inf_{t\in [0,+\infty)}p(t)>0
   \ee
    and
   \be\la{eq1.799}
   I_p\triangleq \int_0^{+\infty}|p'(t)|dt<+\infty.
   \ee
 \begin{theorem}\la{thm1.1} Under the assumptions of Theorem \ref{thm1.2}, suppose that
$p(t)\in C^1[0,\infty)$ satisfies \eqref{eq1.9} and \eqref{eq1.799} in addition.
Then  the initial-boundary-value problem \eqref{1.1}-\eqref{1.8} has a unique global strong solution $(v,u,\te)$ satisfying
 \be
 \begin{cases} \la{1.10}  v ,\,u,\,\theta \in L^\infty(0,\infty;H^1(0,1)),\\ v_t\in
  L^\infty(0,\infty;L^2(0,1))\cap L^2(0,\infty;H^1(0,1)), \\ u_x,\,\theta_x,\,u_t\,,\theta_t,\,v_{xt},\,\,u_{xx},\,\te_{xx} \in
  L^2((0,\infty)\times(0,1)),\end{cases}
  \ee
  and for any $(x,t)\in (0,1)\times (0,+\infty),$
  \bnn
   C^{-1}\leq v(x,t)\leq C ,\quad  C^{-1}\leq \theta(x,t)\leq C
   \enn
   with some positive constant $C$ depending only on the initial data.
  Moreover,  $(v,u,\theta)$ converges to a stationary state $(\hat{v},\hat{u}, 0)$ in $H^1(0,1)$ as $t\rightarrow +\infty$,
  where
  \be\la{eq2.100}
  \hat{v}\triangleq \frac{1}{2\bar{P}}\left\{\int_0^1\left(\frac{u_0^2}{2}+\theta_0+p(0)v_0\right)dx
  +\int_0^{+\infty}p'(\tau)\int_0^1 v(x,\tau)dxd\tau\right\}
  \ee
  and
 \be\la{eq2.10000}
  \hat{\theta}\triangleq \frac{1}{2}\left\{\int_0^1\left(\frac{u_0^2}{2}+\theta_0+p(0)v_0\right)dx
  +\int_0^{+\infty}p'(\tau)\int_0^1 v(x,\tau)dxd\tau\right\}
 \ee
 with \be\la{1.17}
    \bar{P}\triangleq \lim_{t\rightarrow +\infty}p(t).
    \ee
  \end{theorem}

  Similar to the proof of Theorem 3 in \cite{tak}, we have  the following results on  rate of convergence of the strong solutions obtained by Theorem \ref{thm1.1}.

   \begin{theorem}\la{thm1.3} Assume that the assumptions of Theorem \ref{thm1.1} hold. Then there exist positive
  constants $\lambda_1,\lambda_2 (\lambda_1>\lambda_2), C$ depending on $\dis\inf_{t\in [0,+\infty)}p(t), \sup_{t\in [0,+\infty)}p(t)$ and initial data such that the strong solution $(v,u,\te)$  in Theorem \ref{thm1.1} satisfies
  \be\la{eq1.10}
 \int_0^1\{u^2+(v-\hat{v})^2+(\theta-\hat{\theta})^2+v_x^2\}dx\leq C \exp(-\lambda_1 t)
 \left\{\int_0^t \exp(\lambda_1
\tau) F(\tau)d\tau+1\right\},  \ee
 \be\ba\la{eq1.12}
\int_0^1u_x^2dx
&\leq C \exp(-\lambda_1t)\left\{\int_0^t\exp(\lambda_1 \tau)F(\tau)d\tau +1\right\}\\
 &\quad +C\left[\exp(-\lambda_2 t)\left\{\int_0^t \exp(\lambda_2 \tau)(F(\tau)+|p'(\tau)|)d\tau+1\right\}+(p(t)-\bar{P})^2\right],
 \ea\ee
 \be\ba\la{eq1.13}
 \int_0^1 \theta_x^2dx\leq C\exp(-\lambda_2 t)\left\{ \int_0^t \exp(\lambda_2 t)(F(\tau)+|p'(\tau)|)d\tau+1\right\}
 \ea\ee
where
 \be\la{eq1.14}
 Y(t)\triangleq \exp\left\{\int_0^t p(\tau)d\tau\right\},
 \ee
 \bnn
 F(t)\triangleq \frac{1}{Y^2(t)}+\left( \frac{1}{Y(t)}\int_0^t Y(\tau)(\bar{P}-p(\tau))d\tau\right)^2+\int_t^{+\infty}|p'(\tau)|d\tau.
 \enn

  \end{theorem}

We now  make some comments on the analysis of this paper. The key issue is to get the lower and upper bounds of $v$ and $\theta$ (see \eqref{qp1}, \eqref{qp2} and \eqref{eq2.750}). Compared with   \cite{huang} and \cite{tak}, the main difficulties come from the degeneracy and nonlinearity of the heat conductivity due to $\beta>0$ and the outer pressure boundary conditions \eqref{1.8}. Hence, in order to arrive at \eqref{eq2.750},  some new ideas have to be involved here, which are based on three technical treatments. First, applying the standard energetic estimate \eqref{2.6} and the method in  \cite{huang} yields that the specific volume $v$
 is bounded from above and below (see \eqref{2.13} or \eqref{qp1},\eqref{qp2}). Next, the lower bound of $\theta$ (see \eqref{3.91}) is taken into account, and it suffices to  show that $\|(\theta^{-1/2}-\alpha_1^{-1/2})_+\|_{L^\infty[0,1]}\leq C$. To this end, we  multiply \eqref{eq2.4} by $ \tep^p\te^{-3/2}$ and apply Gronwall's inequality  due to the observation that $\dis \int_0^T (\theta^{-1/2}-\al_1^{-1/2})_+^2dt$ can be bounded by $\int_0^T V(t)dt$ (see \eqref{gj4} and \eqref{gj5}).
Finally, for $\beta\in (0,1)$, multiplying \eqref{eq2.4} by $ (\te^{(1-\b)/2}-\al_2^{(1-\b)/2})_+\te^{-(1+\b)/2}$, and together with the standard energetic estimate \eqref{2.6}, we obtain  for any $\beta>0$, \bnn
 \int_0^T \int_0^1 \theta^{-1}\theta_x^2dxdt\leq C,
 \enn
  which not only implies that
$$\int_0^T \max_{x\in [0,1]}\left|\theta(x,t)-\int_0^1 \theta(x,t)dx\right|^2dt\leq C,$$
but also yields that the $L^2((0,T)\times (0,1))$-norm of $\theta_x^2$ is bounded provided $\theta\leq \alpha_2$.
In the case that $\theta>\alpha_2$, we  multiply \eqref{eq2.4} by $(\theta-\al_2)_+$ so that the $L^2((0,T)\times (0,1))$-norm of $\theta_x$ can be bounded by the $L^4(0,T;L^2(0,1))$-norm of $u_x$ which plays an important role in obtaining the
 bound on $L^2((0,T)\times (0,1))$-norm of both $\theta_x$ and $u_{xx}$ (see Lemma \ref{lemma2.9}).  All these discussions  will be carried out in the next section.

\section{Proof of \thmref{thm1.2}}\la{sec2}
We first state   the following   local existence result which can be proved by applying  the contraction mapping principle (c.f. \cite{10,13,tan}).

\begin{lemma}\la{lemma1.1} Let \eqref{1.1}-\eqref{1.8} hold. Then there exists some $T>0$ such that  the initial-boundary-value problem \eqref{1.1}-\eqref{1.8} has a unique strong solution $(v,u,\te)$ satisfying\bnn
 \begin{cases}   v , \,u, \, \theta  \in L^\infty(0,T;H^1(0,1)),\\ v_t\in
  L^\infty(0,T;L^2(0,1))\cap  L^2(0,T;H^1(0,1)), \\ u_t,\,\theta_t,\,v_{xt},\,u_{xx},\,\theta_{xx} \,\in
  L^2((0,T)\times(0,1)).\end{cases}\enn \end{lemma}

  Then, when  $p(t)$ satisfies  \eqref{eq1.9} and \eqref{eq1.799} in addition, we obtain some priori estimates (see \eqref{2.13}, \eqref{2.23},
\eqref{3.31}, \eqref{eq2.750}, \eqref{eq2.86} below), where  the constants depend only on the initial data of the problem so that we  extend the local solution to the whole interval $[0,+\infty)$ and establish Theorem \ref{thm1.1}.

Without loss of generality, we assume that $\ti\mu= \ti\ka=R=c_v=1 $. 
Moreover, integrating \eqref{1.2} over $[0,1]$ and using \eqref{1.8}, then
\be\ba\la{eq200}\int_0^1 udx=\int_0^1 u_0dx.\ea\ee Hence,  we may assume that
\be\ba\la{eq201}\int_0^1 u_0 dx=0.\ea\ee

 Unless otherwise stated, we follow the convention that $C$ is an unspecified positive constant that may vary from expression to expression, even across an inequality (but
not across an equality). Also $C$ depends on initial data, $\beta$ and $p(t)$ generally, and the dependence of $C$ on other parameters will be specified within parenthesis when necessary.

 First, we derive the following representation of $v$ which is essential to arrive at  the  lower  and upper bounds of $v$.

\begin{lemma}\la{lemma3.4}It holds that
\be\ba\la{eq2.7} v(x,t)=\frac{v_0(x)}{B(x,t)Y(t)}\left(1+\frac{1}{v_0(x)}\int_0^t B(x,\tau)\theta(x,\tau)Y(\tau)d\tau\right),
 \ea\ee
where $Y(t)$ is defined in \eqref{eq1.14} and
\bnn
 B(x,t)\triangleq \exp\left\{\int_0^x (u_0(y)-u(y,t))dy\right\}.  \enn
\end{lemma}
\pf By \eqref{1.2}, we get
 \bnn
  u_t=\left(\frac{u_x }{v}-\frac{\theta}{v}\right)_x.
   \enn
Integrating this over $(0,x)$ and using \eqref{1.8} give
\be \la{1.w}
\left(\int_0^xudy\right)_t=\frac{u_x}{v}-\frac{\theta}{v}+p(t).
\ee
Then by \eqref{1.1},
\be\ba\la{eq2.9}
\left(\int_0^xudy\right)_t=\left(\ln v\right)_t-\frac{\te}{v}+p(t).
\ea\ee
Integrating $\eqref{eq2.9}$ over $(0,t)$  yields
\bnn
\ln v=\ln v_0 +\int_0^x(u-u_0)dy+\int_0^t \frac{\theta}{v}d\tau-\int_0^tp(\tau)d\tau,
\enn
 which implies
\bnn\ba
v(x,t)=\frac{v_0(x)}{B(x,t)Y(t)}\exp{\int_0^t \frac{\theta}{v}(x,\tau)d\tau}.
\ea\enn
This in particular gives \eqref{eq2.7} and   finishes the proof of Lemma \ref{lemma3.4}.
\thatsall

Next, we give the following    energetic  estimates.
 \begin{lemma}\la{lemma3.3}  Setting \bnn
V(t)\triangleq \int_0^1 \left(\frac{\theta^\beta \theta_x^2}{v\theta^2}+\frac{u_x^2}{v\theta}\right)(x,t)dx,
\enn we have
\be\ba\la{2.6}
 \sup_{0\le t\le T}\int_0^1\left(u^2/2+v-\ln v+\theta-\ln\theta\right)dx+\int_0^T V(s)ds \leq C_0,
 \ea\ee
 where and in what follows, $C_0$ and $C$ are both   positive constants depending only on $\b, $ $ \|(v_0,u_0,\theta_0)\|_{H^1(0,1)},$ $
 \inf\limits_{x\in [0,1]}v_0(x), $  $ \inf\limits_{x\in [0,1]}\theta_0(x),$ $m_{p,T},$  $M_{p,T},$ and $I_{p,T}$ with \bnn m_{p,T}\triangleq\inf_{t\in [0,T]}p(t),\quad   M_{p,T} \triangleq\sup_{t\in [0,T]}p(t),\quad I_{p,T}\triangleq \int_0^T |p'(t)|dt.\enn

\end{lemma}
\pf
First, integrating \eqref{1.3} over $(0,1)$ and using \eqref{1.1}, \eqref{1.8} immediately  give
\be\ba\la{2.81}
& \left(\int_0^1 \left(\theta +\frac{u^2}{2}\right)dx\right)_t=\left(\frac{ \theta^\beta \theta_x}{v}+\frac{u(u_x-\theta)}{v}\right)\Big{|}_{x=0}^{x=1}\\&=
-p(t)(u(1)-u(0))=-p(t)\int_0^1 u_x dx
\\&=-\left(\int_0^1 p(t) vdx\right)_t+\int_0^1 p'(t)v dx,
 \ea\ee which implies
\be\ba\la{2.5}
\int_0^1 \left(\theta +\frac{u^2}{2}+p(t)v\right)dx
&=\int_0^1 \left(\theta_0 +\frac{u_0^2}{2}+p(0)v_0\right)dx
+\int_0^tp'(\tau)\int_0^1 vdxd\tau\\
&\leq C+\int_0^t |(\ln p(\tau))'| \left(\int_0^1p(\tau) v dx\right) d\tau.
 \ea\ee
Then by  Gronwall's inequality, we have
\be\ba\la{2.3}
\sup_{0\le t\le T}\int_0^1\left(u^2+\theta+v\right)dx\leq C.
 \ea\ee

 Next, by \eqref{1.2}, we rewrite  \eqref{1.3}  as
\be\ba\la{eq2.4}
\theta_t+\frac{\theta}{v}u_x=\left(\frac{\theta^\beta \theta_x}{v}\right)_x+\frac{u_x^2}{v}.
 \ea\ee
Multiplying \eqref{1.1} by $(1-v^{-1})$, \eqref{1.2} by $u$, \eqref{eq2.4} by $1-\theta^{-1}$, and adding them together, we obtain
\be\ba\la{2.8}
& \left(\frac{u^2}{2}+(v-\ln v)+(\theta-\ln\theta)\right)_t+\frac{u_x^2}{v\theta}+\frac{\theta^\beta\theta_x^2}{v\theta^2}
\\&=
\left(u\left(\frac{u_x}{v}-\frac{\theta}{v}\right)\right)_x+u_x+\left((1-\theta^{-1})\frac{\theta^\beta\theta_x}{v}\right)_x.
 \ea\ee
Finally, integrating \eqref{2.8} over $(0,1),$ and  utilizing \eqref{1.1}, \eqref{1.8} lead to
\bnn\ba & \left(\int_0^1(u^2/2+v-\ln v+\theta-\ln\theta)dx\right)_t+\int_0^1\left(\frac{\theta^\beta\theta_x^2}{v\theta^2}+\frac{u_x^2}{v\theta}\right)dx
\\&=
\left((1-p(t))\int_0^1 v dx\right)_t+\int_0^1 p'(t)vdx,
 \ea\enn
which together with \eqref{2.3} yields \eqref{2.6}.
\thatsall

\begin{lemma}\la{lemma3.5}  For any $T>0$, it holds that
\be\ba\la{2.13}
C^{-1}e^{-M_{p,T}T}\leq v(x,t)\leq C,
  \ea\ee for any  $(x,t)\in [0,1]\times [0,T].$
\end{lemma}
\pf
First,  denoting the average of a function $f$ over $(0,1)$ by $$\bar f\triangleq \int_0^1fdx.$$
Since the function $x-\ln x$ is convex, Jensen's inequality gives
$$
\bar \theta -\ln \bar\theta \leq \int_0^1(\theta-\ln \theta)dx,
$$
which together with   \eqref{2.6} implies
\be\ba\la{eq2.13}
\bar{\theta}(t) \in [\alpha_1,\alpha_2],
  \ea\ee
where $0<\alpha_1<\alpha_2$ are two roots of
$$
x-\ln x=C_0,
$$
and $C_0$ is given by Lemma \ref{lemma3.3}.

Thus, combining  \eqref{eq2.13} with \eqref{2.3}, we obtain
\bnn\ba
 \left|\theta^{\frac{\beta+1}{2}}-\bar{\theta}^{\frac{\beta+1}{2}}\right|
 &\leq
\frac{\beta+1}{2}\left(\int_0^1\frac{\theta^\beta \theta_x^2}{\theta^2v}dx\right)^{\frac12} \left(\int_0^1 \theta vdx\right)^{\frac12}\\
& \leq \dis C V^{\frac12}(t)\max_{x\in [0,1]}v^{\frac12}(x,t),
 \ea\enn
  which leads to
\be\ba\la{2.18}
\frac{\alpha_1}{4}-CV(t)\max_{x\in [0,1]}v(x,t)\leq \theta(x,t)\leq C+CV(t)\max_{x\in [0,1]}v(x,t),
\ea\ee
for all $(x,t)\in [0,1]\times [0,T].$

Next, by   \eqref{2.6} and Cauchy's inequality, we have
\bnn
\left|\int_0^xudy\right| \leq\int_0^1|u|dy
\leq\left(\int_0^1u^2dy\right)^{\frac{1}{2}}  \leq C,
\enn
 which  shows that
\be\ba\la{2.15}
 C^{-1}\le  B(x,t)\le  C.
 \ea\ee
Hence, we deduce from \eqref{2.15}, \eqref{eq2.7} and \eqref{2.18} that
\bnn\ba v(x,t)&\leq C+C\int_0^t e^{-\int_\tau^t p(s)ds}\max_{x\in [0,1]}\theta(x,\tau)d\tau
\\&\leq  C+C\int_0^t
 e^{-\int_\tau^t p(s)ds}\left(1+V(\tau)\max_{x\in [0,1]}v(x,\tau)\right)d\tau
\\&
\leq  C+C\int_0^t
 e^{-m_{p,T}(t-\tau)}\left(1+V(\tau)\max_{x\in [0,1]}v(x,\tau)\right)d\tau
\\
&\leq C+C\int_0^t V(\tau) \max_{x\in [0,1]}v(x,\tau)d\tau,
 \ea\enn
which together  with Gronwall's  inequality indicates that
\be\ba\la{2.20}
v(x,t)\leq C,
\ea\ee
for all $(x,t)\in [0,1]\times [0,T].$

Finally, it follows from \eqref{eq2.7} and   \eqref{2.15}  that \be\la{qp3} v(x,t)\ge Ce^{-M_{p,T} T},\ee for all $(x,t)\in [0,1]\times [0,T].$ Combining \eqref{2.20} and \eqref{qp3}, we finish the proof of Lemma \ref{lemma3.5}.
\thatsall

Moreover, if  \eqref{eq1.9} and \eqref{eq1.799} hold,  then all the constants $C$ and $C_0$ stated in Lemmas \ref{lemma3.3} and \ref{lemma3.5} and even the lower bound of  $v$ are indeed independent of $T$. More precisely, we have the following further result.
 \begin{corollary}\la{qwa.1}
If \eqref{eq1.9} and \eqref{eq1.799} hold, then  there exists some positive constant $\ti C$ depending only on $\b, $ $ \|(v_0,u_0,\theta_0)\|_{H^1(0,1)},$ $
 \inf\limits_{x\in [0,1]}v_0(x), $  $ \inf\limits_{x\in [0,1]}\theta_0(x),$ $m_{p},$  $M_{p},$ and $I_p$ such that  \be \la{qp1}  M_{v,T} \triangleq\sup_{t\in [0,T]}\max_{x\in [0,1]} v(x,t)\le \ti C,\ee and
 \be \la{qp2}
  m_{v,T} \triangleq\inf_{t\in [0,T]}\min_{x\in [0,1]} v(x,t)\ge \ti C^{-1}.
 \ee
 \end{corollary}

\pf
If \eqref{eq1.9} and \eqref{eq1.799} hold, then it is easy to see that there exist  positive constants $M_p$ and $\bar{P}$ such that
\be\ba\la{eq7.1}\bar{P}\triangleq
\lim_{t\rightarrow +\infty}p(t)
 \ea\ee
 and
 \be\ba\la{eq7.2} M_p\triangleq\sup_{T\in [0,+\infty)}M_{p,T}=
\sup_{t\in [0,+\infty)}p(t)<\infty.
 \ea\ee
Thus, under such additional condition \eqref{eq1.9} and \eqref{eq1.799}, an analysis completely parallel to those of Lemmas \ref{lemma3.3} and \ref{lemma3.5} except for some necessary modifications when it comes to $p(t)$ shows that all the constants $C_0$ and $C$ in the proof of the above two lemmas are independent of $T$.

Next, it is clear that the inequlity \eqref{qp1} is a direct consequence of \eqref{2.20}, \eqref{eq1.9}, \eqref{eq1.799}, and \eqref{eq7.2}. On the other hand,  it follows from  \eqref{eq2.7},  \eqref{2.6}, \eqref{2.18} and \eqref{eq7.2}  that
\bnn\ba
 v(x,t)&\geq  C\int_0^t e^{-\int_\tau^t p(s)ds}\min_{x\in [0,1]}\theta(x,\tau)d\tau
\\&
\geq C\int_0^t e^{-M_p(t-\tau)}\left(\frac{\alpha_1}{4}-CV(\tau)\right)d\tau
\\&= \frac{C\alpha_1 }{4M_p}-\frac{C\alpha_1 }{4M_p}e^{-M_p t}-C\int_0^t e^{-M_p(t-\tau)}V(\tau)d\tau.
 \ea\enn
Notice that, by \eqref{2.6},
\bnn\ba
\int_0^t e^{-M_p(t-\tau)}V(\tau)d\tau
&=\int_0^{\frac{t}{2}}e^{-M_p(t-\tau)}V(\tau)d\tau+\int_{\frac{t}{2}}^{t}e^{-M_p(t-\tau)}V(\tau)d\tau
\\&
\leq e^{ \frac{-M_pt}{2}}\int_0^\infty V(\tau)d\tau+\int_{\frac{t}{2}}^{t}V(\tau)d\tau
\rightarrow 0, \quad \text{as}\quad t\rightarrow \infty,
 \ea\enn
as a result,  there exists some $\tilde{T}>0$ such that
\bnn
v(x,t)\geq \frac{C\alpha_1}{8M_p},
\enn
for all $(x,t)\in [0,1]\times [\tilde{T}, +\infty).$  Combining this and  \eqref{2.20}, \eqref{qp3}  gives  \eqref{qp1} and \eqref{qp2}. \thatsall

Now, we are in a position to show the lower  bound  of the temperature which is important for our further discussion.


\begin{lemma}\la{lemma3.6}   It holds that for $T>0,$
\be\ba\la{3.91}\inf_{x\in [0,1]}\te(x,t)
\geq C^{-1},
  \ea\ee where $C$ is a  positive constant  depending only on $\b, $ $ \|(v_0,u_0,\theta_0)\|_{H^1(0,1)},$ $
 \inf\limits_{x\in [0,1]}v_0(x), $  $ \inf\limits_{x\in [0,1]}\theta_0(x),$ $m_{p,T},$  $M_{p,T},$  $I_{p,T},$  $m_{v,T}$ and $M_{v,T}. $
\end{lemma}

\begin{remark} \la{rem.1}  It should be note that in the case that $p(t)\in C^1[0,+\infty)$, all the constants $C$ usually depend  on  $\b, $ $ \|(v_0,u_0,\theta_0)\|_{H^1(0,1)},$ $
 \inf\limits_{x\in [0,1]}v_0(x), $  $ \inf\limits_{x\in [0,1]}\theta_0(x),$ $m_{p,T},$  $M_{p,T},$  $I_{p,T},$ and $T.$  However, under additional condition \eqref{eq1.9} and \eqref{eq1.799}, Our discussion below still indicates that $C$ is independent of $T$, which implies that $C$ depends only on initial data,  $\b$, $m_p$, $M_p,$ and $I_p$. In other words, all the estimates in this section will be uniform with respect to time provided \eqref{eq1.9} and \eqref{eq1.799} holds.  \end{remark}

{\it Proof of Lemma \ref{lemma3.6}.}
 For $p>2$, multiplying \eqref{eq2.4} by $ \tep^p\te^{-3/2}$ and integrating  by parts give
\be\ba\la{eq2.26}
&\frac{2}{p+1}\left(\int_0^1 \tep^{p+1}dx\right)_t+ \int_0^1 \frac{u_x^2\tep^p}{v\theta^{3/2}}dx\\
 &\le    \int_0^1 \frac{ u_x  \tep^p}{v\theta^{1/2}}dx\\&\le \frac12 \int_0^1 \frac{u_x^2\tep^p}{v\theta^{3/2}}dx+C  \int_0^1  \theta^{1/2}   \tep^p dx    \\
 &\le \frac12 \int_0^1 \frac{u_x^2\tep^p}{v\theta^{3/2}}dx\\
& \quad+C\max_{x\in [0,1]}(\al_1^{1/2}-\te^{1/2})_+^2\left(1+\int_0^1     \tep^{p+1}  dx \right),
\ea\ee
where we have used the following fact
\bnn\ba & \int_0^1  \theta^{1/2}   \tep^p dx   \\&=  \al_1^{-1}\int_0^1  \theta^{-1/2}   \tep^{p-2}(\al_1^{1/2}-\te^{1/2})_+^2 dx \\&\le C\max_{x\in [0,1]}(\al_1^{1/2}-\te^{1/2})_+^2\left(1+\int_0^1     \tep^{p+1}  dx \right).
\ea\enn

On the other hand,  it follows from \eqref{eq2.13} that for any $\beta\in (0,1)$,
\be\ba\la{gj4}
\max_{x\in [0,1]}(\al_1^{1/2}-\te^{1/2})_+
& \leq C\int_0^1 \theta^{-\frac12}|\theta_x|dx \\
&\leq C \left(\int_0^1 \theta^{\beta-2}\theta_x^2dx\right)^{1/2}\left(\int_0^1 \theta^{1-\beta}dx\right)^{1/2}
 \\
&\leq C V^{1/2}(t)
\ea\ee
and that for $\beta\geq 1$,
\be\ba\la{gj5}
\max_{x\in [0,1]}(\al_1^{1/2}-\te^{1/2})_+
& \leq C\max_{x\in [0,1]}(\al_1^{\b/2}-\te^{\b/2})_+  \\
&\leq C \int_0^1 \theta^{\frac{\beta}{2}-1}|\theta_x|dx \\
&\leq C \left(\int_0^1 \frac{\theta^{\beta}\theta_x^2}{\theta^2 v}dx\right)^{1/2} \\
&\leq C V^{1/2}(t).
\ea\ee
Therefore, by \eqref{eq2.26}-\eqref{gj5} and  Gronwall's inequality, we obtain
\bnn
\|\theta^{-1/2}-\al_1^{-1/2}\|_{L^{p+1}[0,1]}\leq  C,
\enn with $C$ independent of $p.$
Letting $p\rightarrow +\infty$, we get \eqref{3.91}.
\thatsall

Using Lemma \ref{lemma3.6}, we have the following estimate on the $L^\infty(0,T;L^2(0,1))$-norm of $v_x$.
\begin{lemma}\la{lemma3.7} There exists a positive constant $C$ such that
\be\la{2.23}
\sup_{0\leq t\leq T}\int_0^1 v_x^2 dx +\int_0^T \int_0^1 v_x^2(\theta+1)dxdt\leq C.\ee
for any $T>0$.
\end{lemma}

\pf First,
  for $\b\in (0,1)$, multiplying \eqref{eq2.4} by $ (\te^{(1-\b)/2}-\al_2^{(1-\b)/2})_+\te^{-(1+\b)/2}$ gives
\be\ba\la{2.33} &\frac{1}{1-\b}\left(\int_0^1 (\te^{(1-\b)/2}-\al_2^{(1-\b)/2})_+^2dx \right)_t \\
&\quad+\b \int_0^1 \frac{\te_x^2}{ v\te}1_{(\te>\al_2)}dx+ \int_0^1 \frac{u_x^2(\te^{(1-\b)/2}-\al_2^{(1-\b)/2})_+}{v\te^{ (1+\b)/2}}dx\\&= \frac{1+\b}{2}\al_2^{\frac{1-\b}{2}} \int_0^11_{(\te>\al_2)}  \frac{\te^{\frac{\b-3}{2}}\te_x^2}{ v }dx+ \int_0^1 \frac{u_x\te^{ (1-\b)/2}(\te^{(1-\b)/2}-\al_2^{(1-\b)/2})_+}{v}dx\\
&\le \frac{\b}{2} \int_0^1 \frac{\te_x^2}{ v\te}1_{(\te>\al_2)} dx+C\int_0^1 \frac{\theta^\beta \theta_x^2}{v \theta^2 }dx
+ \int_0^1 \frac{u_x\te^{ (1-\b)/2}(\te^{(1-\b)/2}-\al_2^{(1-\b)/2})_+}{v}dx\\
&\le \frac{\b}{2} \int_0^1 \frac{\te_x^2}{ v\te}1_{(\te>\al_2)} dx+C\ve\int_0^1 \frac{u_x^2(\te^{(1-\b)/2}-\al_2^{(1-\b)/2})_+}{v\te^{ (1+\b)/2}}dx+C(\ve) V(t),
\ea\ee
where we have utilized the following two estimates
\bnn\ba &\int_0^1 \frac{|u_x|\te^{ (1-\b)/2}(\te^{(1-\b)/2}-\al_2^{(1-\b)/2})_+}{v}dx\\ &\le \ve\int_0^1 \frac{u_x^2}{v\te^\b}1_{(\te>2\al_2)} dx+C(\ve)\max_{x\in[0,1]}(\te^{(1-\b)/2}-\al_2^{(1-\b)/2})_+^2 \int_0^1 \te dx\\
&\quad +C\int_0^1 \frac{|u_x|(\te^{(1-\b)/2}-\al_2^{(1-\b)/2})_+}{v\te^{1/2}}1_{(\te\le 2\al_2)} dx\\ &\le C\ve\int_0^1 \frac{u_x^2(\te^{(1-\b)/2}-\al_2^{(1-\b)/2})_+}{v\te^{ (1+\b)/2}}dx+C(\ve) V(t)
\ea\enn
and
\be\la{pq01}\ba \max_{x\in[0,1]}(\te^{(1-\b)/2}-\al_2^{(1-\b)/2})_+^2&\le C\max_{x\in[0,1]}(\te^{1/2}-\al_2^{1/2})_+^2\\&\le CV(t).\ea\ee
 Therefore, for $\b\in (0,1),$  choosing $\varepsilon$ suitably small and integrating \eqref{2.33} over $[0,T]$,  we have
\be\la{gj6}
\int_0^T\int_0^1 \theta^{-1}\theta_x^2dxdt\leq C,
\ee which still holds for $\b\ge 1$ due to \eqref{2.6} and  \eqref{3.91}.

Moreover, by \eqref{gj6} and \eqref{eq2.13},
\be\la{gj7}
\int_0^T \max_{x\in [0,1]}(\theta(x,t)-\bar{\theta}(t))^2dt\leq C\int_0^T\left(\int_0^1 \theta^{-1}\theta_x^2dx\right) \left(\int_0^1\theta dx\right)dt\leq C.
\ee
Next, by \eqref{1.1},
  one can rewrite \eqref{1.2} as
\be\la{eq2.40}
\left(u-\frac{v_x}{v}\right)_t-\frac{\theta v_x}{v^2}=-\frac{\theta_x}{v}.
\ee
Multiplying \eqref{eq2.40} by $u-\frac{v_x}{v}$ and integrating over $(0,1)$,  along with \eqref{eq2.13}, we obtain
\bnn\ba
&
\left(\frac{1}{2} \int_0^1 \left(u-\frac{v_x}{v}\right)^2dx\right)_t+\int_0^1 \frac{\theta v_x^2}{v^3}dx\\
&=\int_0^1 \frac{u\theta v_x}{v^2}dx-\int_0^1 \frac{u \theta_x}{v}dx+\int_0^1 \frac{\theta_x v_x}{v^2}dx\\
&\leq C\int_0^1 \frac{u^2\theta}{v}dx+\frac14 \int_0^1 \frac{\theta v_x^2}{v^3}dx
+C \int_0^1 \frac{\theta_x^2}{\theta}dx\\
&\leq C\max_{x\in [0,1]}u^2 \int_0^1 \theta dx+\frac14 \int_0^1 \frac{\theta v_x^2}{v^3}dx+C\int_0^1\frac{\theta_x^2}{\theta}dx\\
&\leq C\max_{x\in [0,1]}u^2+\frac14 \int_0^1 \frac{\theta v_x^2}{v^3}dx +C\int_0^1\frac{\theta_x^2}{\theta}dx,
\ea\enn
which directly leads to
\be\ba\la{eq2.42}
&\int_0^1 \left(u-\frac{v_x}{v}\right)^2dx+\int_0^t\int_0^1 \frac{\theta v_x^2}{v^3}dxdt\\
&\leq C+C\int_0^t \max_{x\in [0,1]}u^2dt +C\int_0^t\int_0^1\frac{\theta_x^2}{\theta}dxdt\\
&\leq C+C\int_0^t V(t)dt\\
&\leq C,
\ea\ee
where in the second inequality  we have used \eqref{gj6},
\be\la{eq202}
|u|=|u-\int_0^1 u dx|\leq \int_0^1 |u_x|dx,
\ee
and
\be\ba\la{eq203}
\int_0^1|u_x|dx&\leq \left(\int_0^1 \frac{u_x^2}{v\theta}dx\right)^{1/2}\left(\int_0^1v\theta dx\right)^{1/2}\\
&\leq C\left(\int_0^1 \frac{u_x^2}{v\theta}dx\right)^{1/2}\leq CV^{1/2}(t).
\ea\ee
 Therefore, we have
\be\la{2.35}
\int_0^1 v_x^2dx+\int_0^T\int_0^1 \theta v_x^2dxdt\leq C,
\ee
due to the following simple fact
$$
\sup_{0\leq t\leq T}\int_0^1 v_x^2dx\leq C\int_0^1 u^2dx+
C\sup_{0\leq t\leq T}\int_0^1 \left(u-\frac{v_x}{v}\right)^2dx\leq C.
$$

Finally,  together with \eqref{3.91} and \eqref{2.35}, we derive
\bnn
\int_0^T\int_0^1 v_x^2dxdt\leq C.
\enn
Combining this with \eqref{2.35} gives \eqref{2.23}, and  we finish the proof of Lemma \ref{lemma3.7}.
\thatsall

For further discussion, we need the following estimate on the $L^2((0,T)\times (0,1))$-norm of $\theta_x$.
\begin{lemma}\la{lemma2.8}  There exists a positive constant $C$ such that
\be\ba\la{3.24}\sup_{0\le t\le T}\int_0^1\te^2dx+
\int_0^T\int_0^1 \theta_x^2dxdt\leq C+C\int_0^T\left(\int_0^1 u_x^2dx\right)^2dt,\ea\ee
for any $T>0.$
\end{lemma}
\pf
Multiplying \eqref{eq2.4} by $(1-\al_2/\te)_+$ and integrating by parts give
\bnn\ba
 &\int_0^1\frac{u_x^2(\te-\al_2 )_+}{v\te }dx  -\left(\int_0^1\int_{\al_2}^\te\left(1-\frac{\al_2}{s} \right)_+dsdx\right)_t\\&=\int_0^1\frac{(\theta-\al_2)_+ }{v}u_x dx+\al_2\int_0^1\frac{\theta^\beta \theta^2_x}{v\te^2} 1_{(\te>\al_2)}dx\\
 &\le \ve\int_0^1\frac{ u_x^2}{v} dx+C(\ve)\int_0^1(\theta-\al_2)_+^2dx+C V(t)\\
 &\le \ve\int_0^1\frac{ u_x^2}{v} dx+C(\ve)\max_{x\in[0,1]}(\theta^{1/2}-\al_2^{1/2} )^2_++C V(t)\\&\le \ve\int_0^1\frac{ u_x^2}{v} dx+C(\ve) V(t) ,
 \ea\enn
 where  in the third inequality we have used \eqref{pq01}.  Consequently,  we get
\be\la{3.26}
\int_0^T\int_0^1 u_x^2dxdt\leq C.
\ee
Multiplying \eqref{eq2.4} by $(\theta-\al_2)_+$ and integrating  by parts, we check that
\be\ba\la{eq2.47}
&\left(\int_0^1\frac{(\theta-\al_2)_+^2}{2}dx\right)_t+\int_0^1 \frac{\theta^\beta \theta_x^2}{v}1_{(\te>\al_2)}dx\\&=-\int_0^1 \frac{\theta(\te-\al_2)_+}{v}u_x dx+\int_0^1 \frac{(\te-\al_2)_+ u_x^2}{v}dx\\
&\leq   C\int_0^1 \frac{ (\te-\al_2)^2_+}{v}|u_x| dx+C\int_0^1 \frac{ (\te-\al_2)_+}{v}|u_x| dx\\
&\quad +\ve\max_{x\in[0,1]}(\te-\al_2)^2_+ +C(\ve) \left(\int_0^1u_x^2dx\right)^2\\
&\leq   C \ve\max_{x\in[0,1]}(\te-\al_2)^2_++C \int_0^1  (\te-\al_2)^2_+ dx\int_0^1  u_x^2 dx \\
&\quad +C(\ve)  \int_0^1u_x^2dx +C(\ve) \left(\int_0^1u_x^2dx\right)^2\\
&\leq C\ve \int_0^1 \theta^\beta \theta_x^21_{(\te>\al_2)}dx+C \int_0^1  (\te-\al_2)^2_+ dx\int_0^1  u_x^2 dx \\
&\quad +C(\ve)  \int_0^1u_x^2dx +C(\ve) \left(\int_0^1u_x^2dx\right)^2,
\ea\ee
where we have used the following inequality
\be\la{eq2.488}
\max_{x\in[0,1]}(\theta-\al_2)_+^2\leq C\int_0^1 \theta^\beta \theta_x^21_{(\te>\al_2)}dx.
\ee
Thus,  choosing $\varepsilon$ suitably small and by Gronwall's inequaltiy, we obtain
\be\ba\la{3.29}
&\sup_{t\in [0,T]}\left(\int_0^1\frac{(\theta-\al_2)_+^2}{2}dx\right)+\int_0^T\int_0^1 \frac{\theta^\beta \theta_x^2}{v}1_{(\te>\al_2)}dxdt\\
&\leq C+C\int_0^T \left( \int_0^1 u_x^2dx\right)^2dt,
\ea\ee
which together with \eqref{2.6} and \eqref{3.91}   gives \eqref{3.24} and the proof is complete.
\thatsall

Now, we can derive the following uniform estimate on the $L^2((0,T)\times (0,1))$-norm of $u_t$ and $u_{xx}$.
\begin{lemma}\la{lemma2.9} There exists a positive constant $C$ such that
\be\la{3.31}
 \sup_{0\leq t\leq T}\int_0^1 u_x^2dx+\int_0^T\int_0^1 \left( u_t^2+u_{xx}^2\right)dxdt\leq C,
 \ee
 for any $T>0$.
\end{lemma}
\pf
First,  by \eqref{1.8},
\bnn
\left(u_x-(\theta-p(t) v)\right)(0,t)=0,\quad \left(u_x-(\theta-p(t)v)\right)(1,t)=0.
\enn
  Multiplying  \eqref{1.2} by  $(u_x-(\theta-p(t)v))_x$ and integrating over
$(0,1)$, by \eqref{1.1} and \eqref{1.3}, we find that
\be\ba\la{2.44}
&\left(\frac12\int_0^1 u_x^2dx\right)_t+\left(\frac{p(t)}{2}\int_0^1 \frac{v_x^2}{v}dx\right)_t
+\int_0^1 \left( \frac{u_{xx}^2}{v}+\frac{\theta_x^2}{v}+p(t)\frac{\theta v_x^2}{v^2}\right)dx
\\
&=\left(\int_0^1 (\theta-p(t)v)u_xdx\right)_t-\frac{p(t)}{2}\int_0^1 u_x\left(\frac{v_x}{v}\right)^2dx
+p(t)\int_0^1 u_x^2dx\\
&\quad+\int_0^1 \frac{\theta^\beta\theta_x}{v}u_{xx}dx-\int_0^1\frac{u_x^3}{v}dx+\int_0^1 \frac{\theta}{v}u_x^2dx-p(t)\int_0^1 \frac{\theta_xv_x}{v}dx\\
&\quad +2\int_0^1\frac{u_{xx}\theta_x}{v}dx
+\int_0^1 \frac{u_{xx}u_xv_x-u_{xx}\theta v_x-u_xv_x\theta_x+\theta v_x \theta_x}{v^2}dx\\
&\quad +\frac12 \int_0^1 p'(t)v\left(\frac{v_x}{v}\right)^2dx+\int_0^1 p'(t)u_x vdx \\
&\leq \left(\int_0^1 (\theta-p(t)v)u_xdx\right)_t+\frac12 \int_0^1 \frac{u_{xx}^2}{v}dx +C\int_0^1\left(u_x^2+v_x^2 + \theta_x^2\right)dx\\
&\quad +C\left(\int_0^1 u_x^2dx\right)^2+C|p'(t)|\left(1+  \int_0^1 u_x^2dx\right)+C\int_0^1 \frac{\theta^{2\beta}\theta_x^2}{v}dx,
\ea\ee
where we have used the following estimates
\bnn\ba
&|\int_0^1 u_x\left(\frac{v_x}{v}\right)^2dx|\\
& \leq C \int_0^1 u_x^2v_x^2dx+C\int_0^1 v_x^2dx\\
&\leq C\max_{x\in [0,1]}u_x^2 \int_0^1 v_x^2dx+C\int_0^1 v_x^2dx\\
&\leq \frac{1}{16}\int_0^1 \frac{u_{xx}^2}{v}dx+C\int_0^1 u_x^2dx+C\int_0^1 v_x^2dx,
\ea\enn
$$
|\int_0^1 \frac{\theta^\beta \theta_x}{v}u_{xx}dx |\leq
\frac{1}{16} \int_0^1  \frac{u_{xx}^2}{v}dx+C\int_0^1 \frac{\theta^{2\beta}\theta_x^2}{v}dx,
$$
\bnn\ba
|\int_0^1 \frac{u_x^3}{v}dx |&\leq C \max_{x\in[0,1]}|u_x|\int_0^1 u_x^2dx\\
&\leq C\max_{x\in[0,1]}u_x^2dx+C\left(\int_0^1 u_x^2dx\right)^2\\
&\leq \frac{1}{16}\int_0^1  \frac{u_{xx}^2}{v}dx+C\int_0^1 u_x^2dx
+C\left(\int_0^1 u_x^2dx\right)^2
\ea\enn
and
$$
|\int_0^1 \frac{\theta}{v}u_x^2dx |\leq C \max_{x\in[0,1]}u_x^2
\leq \frac{1}{16}\int_0^1 \frac{u_{xx}^2}{v}dx+C\int_0^1 u_x^2dx,
$$
$$
|p(t)\int_0^1 \frac{\theta_xv_x}{v}dx |\leq C  \int_0^1 \theta_x^2dx+C\int_0^1 v_x^2dx,
$$
$$
|\int_0^1 \frac{u_{xx}\theta_x}{v}dx |\leq \frac{1}{16} \int_0^1  \frac{u_{xx}^2}{v}dx+C\int_0^1 \theta_x^2dx,
$$
\bnn\ba
&\left|\int_0^1 \frac{u_{xx}u_xv_x-u_{xx}\theta v_x-u_xv_x\theta_x+\theta v_x \theta_x}{v^2}dx\right|\\
&\leq C\int_0^1 |u_{xx}u_xv_x-u_{xx}\theta v_x-u_xv_x\theta_x+\theta v_x \theta_x|dx\\
&\leq \frac{1}{32} \int_0^1  \frac{u_{xx}^2}{v}dx+C\int_0^1\theta_x^2dx+C\int_0^1 u_x^2v_x^2dx+C\int_0^1 \theta^2 v_x^2dx\\
&\leq \frac{1}{32} \int_0^1  \frac{u_{xx}^2}{v}dx+C\int_0^1\theta_x^2dx+C\left(\max_{x\in [0,1]}u_x^2+\max_{x\in [0,1]}(\theta-\bar{\theta})^2+1\right)\int_0^1 v_x^2dx\\
&\leq \frac{1}{16}\int_0^1  \frac{u_{xx}^2}{v}dx+C\int_0^1 u_x^2dx+C\int_0^1\theta_x^2dx+C\int_0^1 v_x^2dx,
\ea\enn
$$
|\int_0^1 p'(t)v\left(\frac{v_x}{v}\right)^2dx|\leq C|p'(t)|\int_0^1 v_x^2dx\leq C|p'(t)|,
$$
$$
|\int_0^1 p'(t)u_x vdx|
\leq C|p'(t)|+C|p'(t)|\int_0^1 u_x^2dx
$$
due to \eqref{2.23} and the following fact
 \be\la{eq2.90}
 \max_{x\in[0,1]}u_x^2\leq \varepsilon \int_0^1 u_{xx}^2dx+C(\varepsilon)\int_0^1 u_x^2dx
 \ee
for any $\varepsilon>0$.

Now, it only remains to   estimate the last term $ \int_0^1 \frac{\theta^{2\beta}\theta_x^2}{v} dx$ on the right side of  \eqref{2.44}.
To this end,  multiplying  \eqref{eq2.4} by $(\theta^{\frac{\beta+2}{2}}-\al_2^{\frac{\beta+2}{2}})_+\theta^{\frac{\beta}{2}}$, we arrive at
\bnn\ba
&\frac{1}{\beta+2}\left(\int_0^1 (\theta^{\frac{\beta+2}{2}}-\al_2^{\frac{\beta+2}{2}})_+^2dx\right)_t+
(\beta+1)\int_0^1 \frac{\theta^{2\beta}\theta_x^2}{v}1_{\{\theta>\alpha_2\}}dx\\
&= \int_0^1\frac{u_x^2\theta^{\frac{\beta}{2}}}{v}(\theta^{\frac{\beta+2}{2}}-\al_2^{\frac{\beta+2}{2}})_+dx
-\int_0^1\frac{u_x \theta^{\frac{\beta+2}{2}}}{v}(\theta^{\frac{\beta+2}{2}}-\al_2^{\frac{\beta+2}{2}})_+dx\\
&\quad +\frac{\beta}{2}\al_2^{\frac{\beta+2}{2}}\int_0^1 \frac{\theta^{\frac{3\beta}{2}}\theta_x^2}{v\theta}1_{\{\theta>\alpha_2\}}dx\\
&\leq \int_0^1\frac{u_x^2\theta^{\frac{\beta}{2}}}{v}(\theta^{\frac{\beta+2}{2}}-\al_2^{\frac{\beta+2}{2}})_+dx
-\int_0^1\frac{u_x \theta^{\frac{\beta+2}{2}}}{v}(\theta^{\frac{\beta+2}{2}}-\al_2^{\frac{\beta+2}{2}})_+dx\\
&\quad +\frac{\beta+1}{2}\int_0^1 \frac{\theta^{2\beta}\theta_x^2}{v}1_{\{\theta>\alpha_2\}}dx+CV(t).
\ea\enn
Consequently,
\be\ba\la{2.45}
&\frac{1}{\beta+2}\left(\int_0^1 (\theta^{\frac{\beta+2}{2}}-\al_2^{\frac{\beta+2}{2}})_+^2dx\right)_t+
\frac{\beta+1}{2}\int_0^1 \frac{\theta^{2\beta}\theta_x^2}{v}1_{\{\theta>\alpha_2\}}dx\\
&\leq \int_0^1\frac{u_x^2\theta^{\frac{\beta}{2}}}{v}(\theta^{\frac{\beta+2}{2}}-\al_2^{\frac{\beta+2}{2}})_+dx
-\int_0^1\frac{u_x \theta^{\frac{\beta+2}{2}}}{v}(\theta^{\frac{\beta+2}{2}}-\al_2^{\frac{\beta+2}{2}})_+dx+CV(t)\\
&\triangleq I_1+I_2+CV(t).
\ea\ee
On the other hand, a straightforward calculation indicates that
\be\ba\la{eq2.68}
I_1&\leq \int_0^1 \frac{u_x^2}{v}(\theta^{\beta+1}-\al_2^{\beta+1})_+dx
\\
&\leq C\max_{x\in [0,1]}\left(\theta^{\beta+1}-\al_2^{\beta+1}\right)_+\int_0^1u_x^2dx\\
&\leq C\int_0^1 |\theta^\beta \theta_x|1_{(\te>\al_2)}dx \int_0^1u_x^2dx\\
&\leq \frac{\beta+1}{4} \int_0^1 \frac{\theta^{2\beta}\theta_x^2}{v}1_{(\te>\al_2)}dx+C\left(\int_0^1 u_x^2dx\right)^2
\ea\ee
and
\be\ba\la{2.47}
I_2&\leq \int_0^1 |u_x|(\theta^{\beta+2}-\al_2^{\beta+2})_+ dx\\
&\leq C\max_{x\in [0,1]}(\theta^{\frac{\beta+2}{2}}-\al_2^{\frac{\beta+2}{2}})_+\int_0^1 (\theta^{\frac{\beta+2}{2}}+
\al_2^{\frac{\beta+2}{2}}) 1_{(\te>\al_2)} |u_x|dx
\\
&\leq C\int_0^1 |\theta^{\frac{\beta}{2}}\theta_x| 1_{(\te>\al_2)}dx \int_0^1 (\theta^{\frac{\beta+2}{2}}+
\al_2^{\frac{\beta+2}{2}})  1_{(\te>\al_2)} |u_x|dx
\\
&\leq C\int_0^1 \theta^\beta \theta_x^2  1_{(\te>\al_2)}dx+C \int_0^1 (\theta^{\frac{\beta+2}{2}}-\al_2^{\frac{\beta+2}{2}})_+^2dx \int_0^1 u_x^2dx+C\int_0^1 u_x^2dx.
\ea\ee
Putting  \eqref{eq2.68}-\eqref{2.47} into \eqref{2.45}, together with \eqref{2.6}, \eqref{3.26} and \eqref{3.29}, and by Gronwall's inequality, we have
\be\ba\la{3.37}
&\sup_{0\leq t\leq T}\left(\int_0^1 (\theta^{\frac{\beta+2}{2}}-\al_2^{\frac{\beta+2}{2}})_+^2dx\right)+\int_0^T \int_0^1 \frac{\theta^{2\beta}\theta_x^2}{v}dxdt\\
&\leq C+C\int_0^T\left(\int_0^1 u_x^2dx\right)^2dt.
\ea\ee

Next,  by \eqref{3.24},
\bnn\ba
\left|\int_0^1 (\theta-p(t)v)u_xdx\right|&\leq  \frac18 \int_0^1 u_x^2dx+C\int_0^1 \theta^2dx+C\\
&\leq \frac18 \int_0^1 u_x^2dx+C\int_0^T\left(\int_0^1 u_x^2dx\right)^2dt+C,
\ea\enn
which together with  \eqref{2.44},  \eqref{3.26},   \eqref{3.24},   \eqref{2.23} and Gronwall's inequality shows that
\be\la{eq2.722}
\int_0^1 u_x^2dx+\int_0^T\int_0^1 u_{xx}^2dxdt\leq C.
\ee
Hence,  by \eqref{eq2.722}, \eqref{3.24}, \eqref{eq2.90} and \eqref{3.37},
\be\la{eq8.88}\sup_{0\le t\le T}\int_0^1 \theta^{\beta+2}dx+
\int_0^T\left(\int_0^1(1+\te^{2\b}) \theta_x^2dx+\max_{x\in [0,1]}u_x^2\right)dt\leq C.
\ee

 Finally, we rewrite \eqref{1.2} as
 \bnn
 u_t=\frac{u_{xx}}{v}-\frac{u_x v_x}{v^2}-\frac{\theta_x}{v}+\frac{\theta v_x}{v^2},
 \enn which along with  \eqref{eq2.722} and  \eqref{eq8.88} implies
 \bnn\ba
 \int_0^T\int_0^1 u_t^2dxdt&\leq  C\int_0^T\int_0^1 u_{xx}^2dxdt+C\int_0^T\int_0^1 u_x^2v_x^2dxdt
\\
 &\quad +C\int_0^T\int_0^1 \theta_x^2dxdt+C\int_0^T\int_0^1 \theta^2v_x^2dxdt\\
 &\leq C+C\int_0^T \max_{x\in [0,1]}u_x^2dt+C\int_0^T \left(\max_{x\in [0,1]}\theta^2\int_0^1 v_x^2dx\right)dt\\
 &\leq C+C\int_0^T \max_{x\in [0,1]}|\theta-\bar{\theta}|^2dt+C\int_0^T\int_0^1 v_x^2dxdt\\
 &\leq C.
 \ea\enn
 Combining this and \eqref{eq2.722}, we finish  the proof of Lemma \ref{lemma2.9}.
 \thatsall

 \begin{lemma}\la{lemma2.10}  There exists a positive constant $C$ such that
\be\la{eq2.750}
C^{-1}\leq \theta(x,t) \leq C,
\ee
for any $(x,t)\in [0,1]\times [0,T].$
\end{lemma}

\pf
Multiplying \eqref{eq2.4} by $\theta^\beta \theta_t$ and integrating  over $(0,1)$, by \eqref{eq8.88}, we give
\be\ba\la{eq2.744}
&\int_0^1 \theta^\beta \theta_t^2dx+\frac12 \left(\int_0^1 \frac{(\theta^\beta\theta_x)^2}{v}dx\right)_t\\
&=-\frac12
\int_0^1\frac{(\theta^\beta\theta_x)^2u_x}{v^2}dx-\int_0^1 \frac{\theta^{\beta+1}\theta_t u_x}{v}dx+\int_0^1 \frac{\theta^\beta \theta_t u_x^2}{v}dx\\
&\leq C\max_{x\in[0,1]}|u_x|\int_0^1(\theta^\beta\theta_x)^2dx+\frac12 \int_0^1 \theta^\beta\theta_t^2dx
+C\int_0^1 \theta^{\beta+2}u_x^2dx\\
&\quad +C\int_0^1\theta^\beta u_x^4dx\\
&\leq C\left(\int_0^1 \frac{(\theta^\beta\theta_x)^2}{v}dx\right)^2+C\max_{x\in [0,1]}u_x^2+\frac12 \int_0^1 \theta^\beta \theta_t^2dx\\
&\quad +C\max_{x\in [0,1]}u_x^2\int_0^1 \theta^{\beta+2}dx
+C\max_{x\in [0,1]}u_x^4\int_0^1\theta^\beta dx\\
&\leq C\left(\int_0^1 \frac{(\theta^\beta\theta_x)^2}{v}dx\right)^2+C\max_{x\in [0,1]}u_x^2+C\max_{x\in [0,1]}u_x^4
+\frac12 \int_0^1 \theta^\beta \theta_t^2dx.
\ea\ee
On the other hand, by \eqref{eq2.722} and \eqref{eq8.88},
\bnn\ba
&\int_0^T\max_{x\in [0,1]}u_x^4dt\\
&\leq C\int_0^T\int_0^1 u_x^4dxdt+C\int_0^T\int_0^1 |u_x^3u_{xx}|dxdt\\
&\leq C\int_0^T\max_{x\in [0,1]}u_x^2\left(\int_0^1 u_x^2dx\right)dt+C\int_0^T \max_{x\in [0,1]}u_x^2\left(\int_0^1 |u_x u_{xx}|dx\right)dt\\
&\leq C+\frac12  \int_0^T \max_{x\in [0,1]}u_x^4dt +C\int_0^T \int_0^1(u_x^2+u_{xx}^2)dx dt\\
&\leq C+\frac12  \int_0^T \max_{x\in [0,1]}u_x^4dt,
\ea\enn
which implies that
\be\la{eq2.766}
\int_0^T\max_{x\in [0,1]}u_x^4dt<C.
\ee
Therefore, together with \eqref{eq2.766}, \eqref{eq8.88} and \eqref{eq2.744}, by Gronwall's inequality, one gets
\be\la{eq2.77}
\sup_{0\leq t\leq T}\int_0^1 (\theta^\beta \theta_x)^2dx+\int_0^T\int_0^1 \theta^\beta \theta_t^2dxdt\leq C.
\ee
In  particular,
\be\ba\la{eq2.79}
&\max_{x\in[0,1]}\left|\theta^{\beta+1}-\bar{\theta}^{\beta+1}\right|\leq (\beta+1)\int_0^1 \theta^\beta |\theta_x|dx\\
&\leq C\left(\int_0^1 (\theta^\beta \theta_x)^2dx\right)^{1/2}\\
&\leq C,
\ea\ee
which  along with \eqref{3.91}  gives  \eqref{eq2.750} and we finish the proof of Lemma \ref{lemma2.10}.
\thatsall

Finally, we obtain  the following  estimate on the $L^2((0,1)\times (0,T))$-norm of $\theta_t$ and $\theta_{xx}$.
\begin{lemma}\la{lemma2.12}  There exists a positive constant $C$ such that
\be\la{eq2.86}
\sup_{0\leq t\leq T}\int_0^1 \theta_x^2dx+\int_0^T\int_0^1 (\theta_x^2+\theta_t^2+\theta_{xx}^2)dxdt\leq C.
\ee
for any $T>0$.
\end{lemma}
\pf  By \eqref{eq2.750} and \eqref{eq2.77},
\be\la{eq2.87}
\sup_{0\leq t\leq T}\int_0^1\theta_x^2dx+\int_0^T\int_0^1\theta_t^2dxdt\leq C.
\ee

Rewrite \eqref{eq2.4} as
\bnn
\frac{\theta^\beta \theta_{xx}}{v}=-\frac{\beta \theta^{\beta-1}\theta_x^2}{v}+\frac{\theta^\beta\theta_xv_x}{v^2}
-\frac{u_x^2}{v}+\frac{\theta u_x}{v}+\theta_t.
\enn
Hence,
\be\ba\la{eq2.89}
\int_0^T\int_0^1 \theta_{xx}^2dxdt &\leq C\int_0^T\int_0^1(\theta_x^4+\theta_x^2v_x^2
+u_x^4+u_x^2+\theta_t^2)dxdt\\
&\leq C+C\int_0^T \max_{x\in [0,1]}\theta_x^2dt+C\int_0^T \max_{x\in [0,1]}u_x^4dt\\
&\leq C+C\int_0^T \int_0^1 |\theta_x\theta_{xx}|dxdt\\
&\leq C+C\int_0^T\int_0^1 \theta_x^2dxdt+\frac12 \int_0^T\int_0^1 \theta_{xx}^2dxdt\\
&\leq C+\frac12 \int_0^T\int_0^1 \theta_{xx}^2dxdt,
\ea\ee
 where we have used \eqref{2.23}, \eqref{3.26}, \eqref{eq2.750}, \eqref{eq2.766} and   \eqref{eq2.87}.
 Therefore, together with \eqref{eq2.89} and \eqref{eq2.87}, we give \eqref{eq2.86} and finish the proof of Lemma \ref{lemma2.12}.
\thatsall

 {\bf Proof of \thmref{thm1.2}}
  With Lemmas \ref{lemma1.1}, \ref{lemma3.5}, \ref{lemma3.6}, \ref{lemma2.9}, \ref{lemma2.10} and  \ref{lemma2.12} at hand, one can prove Theorem \ref{thm1.2} by  a standard argument. We omit it.
\thatsall

\section{Stability of the strong solutions }
In this section, we will  concentrate on studying  the stability of the strong solutions when \eqref{eq1.9} and \eqref{eq1.799} both  hold. As indicated in Remark \ref{rem.1}, all the constants $C,$ $\ti C,$ and $C_0$ stated in Section 2 are independent of $T.$

\begin{lemma}\la{lemma3.1}  It holds that
\be\la{3.1}
\|(v-\hat{v}, u,\theta-\hat{\theta})(\cdot,t)\|_{H^1(0,1)}\rightarrow 0, \quad \text{as} \quad  t\rightarrow +\infty.
\ee
where $\hat{v}, \hat{\theta}$ are defined in \eqref{eq2.100},  \eqref{eq2.10000} respectively.
\end{lemma}

\pf First,  we check that
\be\la{eq2.201}
\int_0^{+\infty} \left(\left|\frac{d}{dt}\|v_x(\cdot,t)\|_{L^2}^2\right|
+\left|\frac{d}{dt}\|\theta_x(\cdot,t)\|_{L^2}^2\right|\right)dt
\leq C,
\ee
where we have taken advantage of the following facts
\bnn
\int_0^1 v_x v_{xt}dx=\int_0^1 v_x u_{xx}dx
\enn
and
\bnn
\int_0^1 \theta_x \theta_{xt}dx=-\int_0^1 \theta_{xx} \theta_t dx.
\enn
As a result, we deduce from \eqref{2.23}, \eqref{eq2.86} and \eqref{eq2.201} that
\be\la{eq2.202}
\lim_{t\rightarrow +\infty}\int_0^1( v_x^2+ \theta_x^2)dx=0.
\ee
In particular,
\be\la{3.4}
 \lim_{t\rightarrow +\infty}\int_0^1( (v-\bar v)^2+ (\te-\bar \te)^2)dx\le C\lim_{t\rightarrow +\infty}\int_0^1( v_x^2+ \theta_x^2)dx=0.
\ee

Multiplying \eqref{1.2} by $u$   and integrating  over $[0,1]$
lead to
\bnn
\frac12\left(\int_0^1 u^2dx\right)_t+\int_0^1 \frac{u_x^2}{v}dx=\int_0^1 \left(\frac{\theta}{v}-p(t)\right)u_xdx
\leq C\int_0^1 |u_x|dx,
\enn
 which yields that there exists some constant $C>0$ such that  for any $N>0$
  and $s,t\in [N,N+1]$
  $$
  \int_0^1 u^2(x,t)dx-\int_0^1 u^2(x,s)dx\leq C\int_N^{N+1}\int_0^1 |u_x|dxd\tau.
  $$
  Integrating this with respect to $s$ over $(N,N+1)$, together with \eqref{eq202}, we get
  $$
  \sup_{N\leq t\leq N+1}\int_0^1 u^2(x,t)dx\leq C\left(\int_N^{N+1}\int_0^1 u_x^2dxd\tau\right)^{1/2}.
  $$

  Letting $N\rightarrow +\infty$ and by \eqref{3.26},
  \be\la{3.5}
  \lim_{t\rightarrow +\infty}\int_0^1 u^2(x,t)dx=0.
  \ee
Next, multiplying \eqref{1.w} by $v$ and integrating  over $[0,1]$, we obtain
\be \la{2.w}
\bar\te-p(t)\bar v=  \int_0^1u_x dx+\int_0^1u_x \int_0^xudy dx-\left(\int_0^1v\int_0^xudydx\right)_t,
\ee
then
\be\ba\la{3.6}
\int_0^t(\bar\te-p(t)\bar v)^2dt&\le \frac12\int_0^t(\bar\te-p(t)\bar v)^2dt+C\int_0^t \int_0^1u^2_x dxdt +C \\&\quad +C\int_0^t\int_0^1 | u | dx\left( |\bar\te'|+|p'(t)|\bar v+p(t)|\bar v'| \right) dt.
\ea\ee

On the other hand, Observe that, by \eqref{eq2.4},
\bnn \bar\te'=-\int_0^1\frac{\te}{v}u_xdx+\int_0^1\frac{u_x^2}{v} dx, \enn
which along with  \eqref{3.26}, \eqref{3.31} and \eqref{eq2.750} yields
\be \ba\la{3.8} \int_0^{+\infty} (\bar\te')^2dt\le C\ea \ee
and
\be  \ba\la{3.9} \int_0^{+\infty} (\bar v')^2dt=\int_0^{+\infty} \left(\int_0^1u_{x}dx\right)^2dt\le C. \ea \ee
Therefore,
\be  \ba\la{3.11}
\int_0^{+\infty}(\bar\te-p(t)\bar v)^2dt \le  C. \ea \ee
Moreover, by  \eqref{3.8}- \eqref{3.11},
\be  \ba\la{3.12}
&
\int_0^{+\infty}\left|\frac{d}{dt}(\bar\te-p(t)\bar v)^2\right|dt\\ & \le  C
\int_0^{+\infty}\left|  \bar\te-p(t)\bar v \right|\left( |\bar\te'|+|p'(t)|\bar v+p(t)|\bar v'| \right)dt\\& \le  C
\int_0^{+\infty}\left(  \bar\te-p(t)\bar v \right)^2dt+C\int_0^{+\infty}\left( (\bar\te')^2+|p'(t)|+(\bar v')^2 \right)dt \\&\le C, \ea \ee
which together with \eqref{3.11} leads to
\bnn\lim_{t\rightarrow +\infty} (\bar\te-p(t)\bar v)=0. \enn
In addition, by \eqref{2.5}, we have
\be\ba\la{3.121}
\bar\te+p(t)\bar v =-\frac12\int_0^1u^2dx
+\int_0^1 \left(\theta_0 +\frac{u_0^2}{2}+p(0)v_0\right)dx
+\int_0^tp'(\tau)\int_0^1 vdxd\tau,
 \ea\ee
Hence, by \eqref{3.5} and \eqref{3.121},
\be  \ba \la{3.15}
\lim_{t\rightarrow +\infty}  \bar v  =\hat v,\quad \lim_{t\rightarrow +\infty}   \bar\te  =\hat \te.
 \ea \ee

Finally,  introducing
$$
\Phi(t)\triangleq \int_0^1\left( u_x  -  \te+p(t)v \right)^2dx.
$$
Multiplying \eqref{1.w} by $v$,  one checks that
\bnn
u_x-\theta+p(t)v=v\left(\int_0^x udx\right)_t,
\enn
which together with \eqref{3.31} implies
\be\ba\la{3.16}
\int_0^{+\infty} \Phi(t) dt \leq C\int_0^{+\infty} \int_0^1 u_t^2dxdt\leq C.
\ea\ee
On the other hand,  by \eqref{1.8}, \eqref{3.26}, \eqref{3.31}, \eqref{eq2.86} and \eqref{3.16}, we conclude that
\bnn  \ba \int_0^{+\infty}|\Phi'(t)|dt\le & C  \int_0^{+\infty}\left|\int_0^1\left( u_x-\te+p(t)v \right)  u_{xt}   dx\right|dt\\&+ C  \int_0^{+\infty}\left|\int_0^1\left( u_x  -  \te+p(t)v \right)\left(  -  \te+p(t)v \right)_t dx\right|dt\\= & C  \int_0^{+\infty}\left|\int_0^1\left( u_x  -  \te+p(t)v \right)_x  u_{t}   dx\right|dt\\&+ C  \int_0^{+\infty}\left|\int_0^1\left( u_x  -  \te+p(t)v \right)\left(  -  \te+p(t)v \right)_t dx\right|dt\\
\leq& C\int_0^{+\infty} \int_0^1 u_{xx}^2dxdt +C\int_0^{+\infty} \int_0^1 \theta_x^2dxdt
+C\int_0^{+\infty} \int_0^1 v_x^2dxdt\\
&+C\int_0^{+\infty} \int_0^1 u_t^2dxdt+C\int_0^{+\infty}\int_0^1 u_x^2dxdt\\
&+C\int_0^{+\infty} \int_0^1 \theta_t^2dxdt+C\int_0^{+\infty} |p'(t)|dt\\
\leq &C.
\ea \enn
Combining the above inequality with \eqref{3.16}, we derive
\be\la{3.18}
\lim_{t\rightarrow +\infty }\Phi(t)=0,
\ee
which together with \eqref{eq2.202}  gives
\be\la{3.19}
\lim_{t\rightarrow +\infty}\int_0^1 u_x^2dx=0.
\ee
As a result, by virtue of \eqref{eq2.202},  \eqref{3.4}, \eqref{3.5}, \eqref{3.15} and \eqref{3.19}, we establish  \eqref{3.1}.
\thatsall

{\bf Proof of Theorem \ref{thm1.1}}  Notice that all the conclusions  in  Section \ref{sec2} are uniform with respect to time provided that \eqref{eq1.9} and \eqref{eq1.799} holds,  by Theorem \ref{thm1.2}, the problem \eqref{1.1}-\eqref{1.8} has a unique global strong solution $(v,u,\te)$ satisfying \eqref{1.111} and $v$, $\theta$ has a uniform positive below and upper bounds, Moreover, by Lemma \ref{lemma3.1}, such $(v,u,\theta)$ converges to a stationary state $(\hat{v},\hat{u},0)$ in $H^1(0,1)$ as $t\rightarrow +\infty$. we finish the proof of Theorem \ref{thm1.1}.

 \end{document}